\documentclass[letterpaper, 10pt, conference]{ieeeconf}  

\IEEEoverridecommandlockouts                              
\usepackage[utf8]{inputenc}
\usepackage[T1]{fontenc}
\usepackage{amsmath}
\usepackage{mathrsfs}
\usepackage{amsfonts}
\usepackage{algorithm}
\usepackage{algorithmicx}
\usepackage{algpseudocode}

\usepackage{color}

 \usepackage{graphicx}

\title{\LARGE \bf
Efficient solution method based on inverse dynamics for \\ optimal control problems of rigid body systems
}

\author{Sotaro Katayama$^{1}$ and Toshiyuki Ohtsuka$^{1}$
\thanks{$^{1}$S. Katayama and T. Ohtsuka are with Department of System Science, Graduate School of Informatics, Kyoto University, Kyoto, Japan
        {\tt\small katayama.25w@st.kyoto-u.ac.jp}, 
        {\tt\small ohtsuka@i.kyoto-u.ac.jp}}%
}

\begin{document}


\onecolumn
\noindent
© 2021 IEEE. Personal use of this material is permitted. Permission from IEEE must be obtained for all other uses, in any current or future media, including reprinting/republishing this material for advertising or promotional purposes, creating new collective works, for resale or redistribution to servers or lists, or reuse of any copyrighted component of this work in other works.

\hspace{1cm}

\noindent
\textbf{Published article:} \\ 
\noindent
S. Katayama and T. Ohtsuka, ``Efficient solution method based on inverse dynamics for optimal control problems of rigid body systems,'' 2021 IEEE International Conference on Robotics and Automation (ICRA), 2021, pp. 2070--2076, doi: 10.1109/ICRA48506.2021.9561109.

\twocolumn

\maketitle
\thispagestyle{empty}
\pagestyle{empty}

\begin{abstract}

We propose an efficient way of solving optimal control problems for rigid-body systems on the basis of inverse dynamics and the multiple-shooting method. We treat all variables, including the state, acceleration, and control input torques, as optimization variables and treat the inverse dynamics as an equality constraint. We eliminate the update of the control input torques from the linear equation of Newton's method by applying condensing for inverse dynamics. The size of the resultant linear equation is the same as that of the multiple-shooting method based on forward dynamics except for the variables related to the passive joints and contacts. Compared with the conventional methods based on forward dynamics, the proposed method reduces the computational cost of the dynamics and their sensitivities by utilizing the recursive Newton-Euler algorithm (RNEA) and its partial derivatives. In addition, it increases the sparsity of the Hessian of the Karush–Kuhn–Tucker conditions, which reduces the computational cost, e.g., of Riccati recursion. Numerical experiments show that the proposed method outperforms state-of-the-art implementations of differential dynamic programming based on forward dynamics in terms of computational time and numerical robustness. 

\end{abstract}

\section{Introduction}
Optimal control plays a significant role in motion planning and control such as trajectory optimization (TO) and model predictive control (MPC) \cite{bib:nmpc} for rigid-body systems. TO can generate dynamically consistent motion under versatile control objectives and constraints even for highly nonlinear systems such as underactuated systems by solving the optimal control problem (OCP). MPC leverages the same advantages as TO for real-time control by solving the OCP at each sampling time. However, we still have to improve the computational efficiency and numerical robustness of the OCP for rigid-body systems whose dynamics are complicated and highly nonlinear.

The computational time of the OCP for rigid-body systems depends substantially on the computational time of the dynamics and their sensitivities. Most previous researches on the OCP for rigid-body systems (e.g., \cite{bib:legged, bib:legged2, bib:crocoddyl}) have been based on forward dynamics, which is a calculation of the generalized acceleration for the given configuration, generalized velocity, and generalized torques. These studies incorporate forward dynamics into the state equation, which is a natural representation, especially for single-shooting methods such as differential dynamic programming (DDP) \cite{bib:ddp} and the iterative linear quadratic regulator (iLQR) \cite{bib:ilqr}. They utilize the articulated body algorithm (ABA) \cite{bib:aba}, an efficient recursive algorithm to compute forward dynamics, or directly compute the inverse of the joint inertia matrix, in solving the OCP. In contrast, our previous work \cite{bib:inverseocp} proposed basing the OCP on inverse dynamics to reduce the computational time compared with that of the OCP based on forward dynamics. Inverse dynamics are calculations of the generalized torques for a given configuration, generalized velocity, and generalized acceleration. We can compute the inverse dynamics with a smaller computational time than with forward dynamics because the recursive Newton-Euler algorithm (RNEA) is faster than ABA \cite{bib:featherstone}. 
Moreover, we can compute the sensitivities of inverse dynamics faster than those of forward dynamics \cite{bib:analytical:rbd, bib:numerical:diff}. 
Therefore, we can reduce the computational cost of the OCP for rigid-body systems by using inverse dynamics instead of forward dynamics, as illustrated numerically in \cite{bib:inverseocp}. 
However, our previous approach of \cite{bib:inverseocp} has a potential drawback when it is utilized for MPC because it is a single-shooting method; i.e., it regards only the acceleration as the decision variable, which means it cannot leverage parallel computing and can lack numerical robustness. 

The single-shooting method regards only the control input as decision variables and computes the state on the horizon by simulating the system's dynamics based on the control input at each iteration. Typical examples of the single-shooting method are DDP and iLQR, which solve the linear equation of Newton's method by using Riccati recursion \cite{bib:Riccati} and are very popular in robotics applications \cite{bib:legged, bib:legged2, bib:crocoddyl}. 
In contrast, the multiple-shooting method regards all variables (the state, costate, and control inputs) as optimization variables and can inherently leverage parallel computing by splitting the computation of the residual of the Karush–Kuhn–Tucker (KKT) conditions and the Hessian of the KKT condition into each time stage, which is impossible for the single-shooting methods due to their expensive serial computational parts, such as the simulation of the system's dynamics over the horizon. 
Even in a single-thread computation, the computational cost of the multiple-shooting method is almost the same as that of the single-shooting method even though it has more optimization variables thanks to structure-exploiting Newton-type methods \cite{bib:condensing, bib:liftedNewton, bib:Riccati}. 
Further, the multiple-shooting method is empirically known to converge quickly even if the initial guess of the solution is far from the (local) optimal solution; In contrast, the single-shooting method converges slowly or even diverges in such a situation.

In this paper, we propose an efficient solution method of the OCP for rigid-body systems based on inverse dynamics and the multiple-shooting method. We treat all variables, including the state, acceleration, and control input torques, as optimization variables, and treat the inverse dynamics as an equality constraint. We eliminate the update of the control input torques from the linear equation of Newton's method by applying condensing for inverse dynamics. The size of the resultant linear equation is the same as that of the multiple-shooting method based on forward dynamics except for the variables related to the passive joints and contacts. Compared with the conventional methods based on forward dynamics, the proposed method reduces the computational cost of the dynamics and their sensitivities by utilizing RNEA and its partial derivatives. In addition, it increases the sparsity of the Hessian of the KKT conditions, which reduces the computational cost, e.g., of Riccati recursion.
Note that the inverse dynamics-based formulation is also utilized in a contact-implicit TO \cite{bib:todorov:invdyn, bib:todorov:mujoco} to alleviate the numerical ill-conditioning due to an approximation with smooth contact model. 
In contrast to these studies, our method is not limited to a specific contact-implicit TO, e.g., it can treat rigid contacts, as illustrated in a numerical experiment.

This paper is organized as follows. In Section \ref{section:formulation}, we formulate the OCP based on inverse dynamics. Section \ref{section:solutionmethod} introduces the proposed solution method with condensing of inverse dynamics. Section \ref{section:expriment} compares the proposed method with state-of-the-art implementations of DDP/iLQR and demonstrates its effectiveness in terms of the computational time and numerical robustness. We conclude in Section \ref{section:conclu} with a brief summary and mention of future work. 

\textit{Notation:} We describe the partial derivatives of a differentiable function by certain variables using a function with subscripts; i.e., $f_x (x)$ denotes $\frac{\partial f}{\partial x} (x)$ and $g_{x y} (x, y)$ denotes $\frac{\partial^2 g}{\partial x \partial y} (x, y)$. We denote a diagonal matrix whose elements are a vector $x$ by ${\rm diag} (x)$. Furthermore, we denote an $n \times n$ identity matrix by $I_{n}$ and an $n \times m$ zero matrix by $O_{n \times m}$.

\section{Optimal Control Problem Based on Inverse Dynamics}\label{section:formulation}
\subsection{Rigid-body systems}
Let $Q$ be the configuration manifold of the rigid-body system. Let $q \in Q$, $v \in \mathbb{R}^n$, $a \in \mathbb{R}^n$, $f \in \mathbb{R}^{n_f}$, and $u \in \mathbb{R}^n$ be the configuration, generalized velocity, acceleration, stack of the contact forces, and input torques, respectively. The equation of motion of the rigid-body system is given by
\begin{equation}\label{eq:equationofmotion}
  M(q) a + h(q, v) - J^{\rm T} (q) f = u,
\end{equation}
where $M (q) \in \mathbb{R}^{n \times n}$ denotes the inertia matrix, $h (q, v) \in \mathbb{R}^{n}$ encompasses Coriolis, centrifugal, and gravitational terms, and $J (q) \in \mathbb{R}^{n_f \times n}$ denotes the stack of the contact Jacobians. We write (\ref{eq:equationofmotion}) as an equality constraint of the inverse dynamics as follows:
\begin{equation}\label{eq:id}
    {\rm ID} (q, v, a, f) - u = 0,
\end{equation}
where ${\rm ID} (q, v, a, f)$ is defined by the left-hand side of (\ref{eq:equationofmotion}) in the formulation of the OCP. Note that we can efficiently compute ${\rm ID} (q, v, a, f)$ by using RNEA and its partial derivatives by using the partial derivatives of RNEA \cite{bib:analytical:rbd}. We also assume that the input torques $u$ are an $n$-dimensional vector even when the system is underactuated, as RNEA and its partial derivatives are well-defined only for fully actuated systems. As stated in the next subsection, we treat the underactuated systems by introducing an equality constraint that zeros the elements of $u$ corresponding to passive joints.

\subsection{Optimal control problem}
We consider an OCP with $N$ time stages. 
The configuration, velocity, acceleration, external forces, and input torques for $N$ stages are denoted as $q_0, ..., q_N \in Q$, $v_0, ..., v_N \in \mathbb{R}^n$, $a_0, ..., a_{N-1} \in \mathbb{R}^n$, $f_0, ..., f_{N-1} \in \mathbb{R}^{n_f}$, $u_0, ..., u_{N-1} \in \mathbb{R}^n$, all of which are regarded as optimization variables to formulate the OCP based on inverse dynamics.
We assume the initial state is given by $\bar{q}$ and $\bar{v}$, and then consider the constraints on the initial state  
\begin{equation}\label{eq:initialstate}
    \delta(\bar{q}, q_0) = 0, \; \bar{v} - v_0  = 0, \;\;\; \bar{q} \in Q, \; \bar{v} \in \mathbb{R}^n,
\end{equation}
where $\delta (q_1, q_2) \in \mathbb{R}^n$ denotes the subtraction operation of the configurations $q_2 \in Q$ from $q_{1} \in Q$ on the manifold $Q$. The state equation discretized with the forward Euler method is given by 
\begin{align}\label{eq:forwardEuler}
    \begin{bmatrix}
        \delta (q_i, q_{i+1}) + v_i \Delta \tau \\
        v_i - v_{i+1} + a_i \Delta \tau
    \end{bmatrix} = 0, \;\; i=0, ..., N-1,
\end{align}
where $\Delta \tau$ is the time step of the discretization given by $\Delta \tau = T/N$ with the length of the horizon $T$. The state equation (\ref{eq:forwardEuler}) takes a simple form since we consider the equation of the motion of the system (\ref{eq:equationofmotion}) as an equality constraint (\ref{eq:id}) in the OCP and regard the acceleration $a_i$ as the optimization variables. 
In general, we can also assume an $m_c$-dimensional equality constraint, 
\begin{equation}\label{eq:equalityconstraints}
    C (q_i, v_i, a_i, u_i, f_i) \Delta \tau = 0 , \;\; i=0, ..., N-1,
\end{equation}
and $m_g$-dimensional inequality constraint, 
\begin{equation}\label{eq:inequalityconstraints}
    g (q_i, v_i, a_i, u_i, f_i) \Delta \tau \leq 0, \;\; i=0, ..., N-1.
\end{equation}
Note that we multiply $C (\cdot)$ and $g (\cdot)$ by $\Delta \tau$ in (\ref{eq:equalityconstraints}) and (\ref{eq:inequalityconstraints}) so that the proposed formulation will correspond to the continuous-time Euler-Lagrange equations discretized with a time step $\Delta \tau$ \cite{bib:appliedOCP}.
If the system is underactuated, (\ref{eq:equalityconstraints}) contains the elements of the control input torques corresponding to the passive joints. For example, if the system has a floating base whose joint indices are 1 to 6, $C (q_i ,v_i, a_i, u_i, f_i)$ in (\ref{eq:equalityconstraints}) includes $[u_i ^{(1)} \, u_i ^{(2)} \, u_i ^{(3)} \, u_i ^{(4)} \, u_i ^{(5)} \, u_i ^{(6)} ]^{\rm T}$, where $u_i ^{(j)}$ is the $j$-th element of $u_i$. 

For the above description of the rigid-body system and constraints, the OCP is given by 
\begin{equation*}\label{eq:cost}
    \min_{q_i, v_i, a_i, u_i, f_i} J = \varphi(q_N, v_N) + \sum_{i=0}^{N-1} l(q_i, v_i, a_i, u_i, f_i) \Delta \tau, 
\end{equation*}
subject to (\ref{eq:id})--(\ref{eq:inequalityconstraints}), where $\varphi(q_N, v_N)$ denotes the terminal cost and $l(q_i, v_i, a_i, u_i, f_i) \Delta \tau$ denotes the stage cost. 
We treat the inequality constraints (\ref{eq:inequalityconstraints}) by using the primal-dual interior point method (PDIPM) \cite{bib:nocedal}, which can treat large-scale and nonlinear constraints efficiently.  
The inequality constraint (\ref{eq:inequalityconstraints}) is then transformed into an equality constraint by introducing slack variables $s_i \in \mathbb{R}^{m_g}$,
\begin{equation}\label{eq:slack}
    g (q_i, v_i, a_i, u_i, f_i) \Delta \tau + s_i \Delta \tau = 0,
\end{equation}
with an additional inequality constraint $s_i \geq 0$. 

\subsection{KKT conditions}
Next, we derive the KKT conditions, necessary conditions for optimal control \cite{bib:appliedOCP, bib:nocedal}.
Let the set of primal variables at stage $i$ be $y_i := \left\{ q_i \;\; v_i \;\; a_i \;\; f_i \;\; u_i \right\}$ and the set of primal variables without the control input torques be $\tilde{y}_i := \left\{ q_i \;\; v_i \;\; a_i \;\; f_i \right\}$ for $i = 0, ..., N-1$. We then introduce the Lagrangian of this OCP with PDIPM,
\begin{align*}\label{eq:Lagrangian}
    & \mathcal{L} = \varphi(q_N, v_N) + \sum_{i=0}^{N-1} l(y_i) \Delta \tau +  \begin{bmatrix} \lambda_{0} \\ \gamma_{0} \end{bmatrix} ^{\rm T} \begin{bmatrix} \delta(\bar{q}, q_0) \\ \bar{v} - v_0 \end{bmatrix} \notag \\ 
    & + \sum_{i=0}^{N-1} \begin{bmatrix} \lambda_{i+1} \\ \gamma_{i+1} \end{bmatrix}  ^{\rm T} 
    \begin{bmatrix} \delta(q_i, q_{i+1}) + v_i \Delta \tau \\ v_i - v_{i+1} + a_i \Delta \tau \end{bmatrix} \notag \\
    & + \sum_{i=0}^{N-1} \beta_{i} ^{\rm T} ({\rm ID} (\tilde{y}_i) - u_i) \Delta \tau + \sum_{i=0}^{N-1} \mu_i ^{\rm T} C(y_i) \Delta \tau \notag \\ 
    & + \sum_{i=0}^{N-1} \nu_i ^{\rm T} (g(y_i) \Delta \tau + s_i \Delta \tau ) - \sum_{i=0}^{N-1} \epsilon \ln s_i \Delta \tau, 
\end{align*}
where $\lambda_i, \; \gamma_i \in \mathbb{R}^{n}$, $\mu_i \in \mathbb{R}^{m_c}$, and $\nu_i \in \mathbb{R}^{m_g}$ denote the Lagrange multipliers with respect to the state equation (\ref{eq:forwardEuler}), the equality constraint (\ref{eq:equalityconstraints}), and the inequality constraints (\ref{eq:inequalityconstraints}), and $\epsilon > 0$ denotes the barrier parameter. The KKT conditions are then given by (\ref{eq:id})--(\ref{eq:equalityconstraints}), (\ref{eq:slack}),
\begin{equation}\label{eq:LqvN}
    \begin{bmatrix}
        \mathcal{L}_{q_N} ^{\rm T} \\
        \mathcal{L}_{v_N} ^{\rm T} 
    \end{bmatrix}
    = \begin{bmatrix}
        \varphi_{q_N} ^{\rm T} (q_N, v_N) \\
        \varphi_{v_N} ^{\rm T} (q_N, v_N)
    \end{bmatrix}
    + \begin{bmatrix}
        \delta_{q_{N}} ^{\rm T} (q_{N-1}, q_{N}) \lambda_{N} \\
        - \gamma_{N}
    \end{bmatrix} = 0,
\end{equation}
and the following equations for $i=0, ..., N-1$,
\begin{align}\label{eq:Lqvi}
    \begin{bmatrix}
        \mathcal{L}_{q_i} ^{\rm T} \\
        \mathcal{L}_{v_i} ^{\rm T} 
    \end{bmatrix}
    = & \;
    \begin{bmatrix}
        l_{q_i} ^{\rm T} (y_i) \Delta \tau \\
        l_{v_i} ^{\rm T} (y_i) \Delta \tau 
    \end{bmatrix}
    + \begin{bmatrix}
        \delta_{q_{i}} ^{\rm T} (q_{i}, q_{i+1}) & O_{n \times n} \\
        I_n \Delta \tau & I_n 
    \end{bmatrix}
    \begin{bmatrix}
        \lambda_{i+1} \\
        \gamma_{i+1}
    \end{bmatrix} \notag \\
    & + \begin{bmatrix}
        {\rm ID}_{q_i} ^{\rm T} \Delta \tau (\tilde{y}_i) \\
        {\rm ID}_{v_i} ^{\rm T} \Delta \tau (\tilde{y}_i)
    \end{bmatrix}
    \beta_i  
    + \begin{bmatrix}
        {C}_{q_i} ^{\rm T} \Delta \tau ({y}_i) \\
        {C}_{v_i} ^{\rm T} \Delta \tau ({y}_i)
    \end{bmatrix}
    \mu_i \notag \\
    & + \begin{bmatrix}
        {g}_{q_i} ^{\rm T} \Delta \tau ({y}_i) \\
        {g}_{v_i} ^{\rm T} \Delta \tau ({y}_i)
    \end{bmatrix}
    \nu_i 
    + \begin{bmatrix}
        \delta_{q_{i}} ^{\rm T} (q_{i-1}, q_i) \lambda_{i} \\
        - \gamma_{i}
    \end{bmatrix} = 0,
\end{align}
\begin{align}\label{eq:Lafi}
    \begin{bmatrix}
        \mathcal{L}_{a_i} ^{\rm T} \\
        \mathcal{L}_{f_i} ^{\rm T} 
    \end{bmatrix}
    = & \;
    \begin{bmatrix}
        l_{a_i} ^{\rm T} (y_i) \Delta \tau \\
        l_{f_i} ^{\rm T} (y_i) \Delta \tau 
    \end{bmatrix}
    + \begin{bmatrix}
        O_{n \times 1} \\
        \gamma_{i+1} \Delta \tau
    \end{bmatrix}
    + \begin{bmatrix}
        {\rm ID}_{a_i} ^{\rm T} (\tilde{y}_i) \Delta \tau \\
        {\rm ID}_{f_i} ^{\rm T} (\tilde{y}_i) \Delta \tau 
    \end{bmatrix}
    \beta_i 
    \notag \\
    & + \begin{bmatrix}
        {C}_{a_i} ^{\rm T} ({y}_i) \Delta \tau \\
        {C}_{f_i} ^{\rm T} ({y}_i) \Delta \tau 
    \end{bmatrix}
    \mu_i
    + \begin{bmatrix}
        {g}_{a_i} ^{\rm T} ({y}_i) \Delta \tau \\
        {g}_{f_i} ^{\rm T} ({y}_i) \Delta \tau
    \end{bmatrix}
    \nu_i = 0,
\end{align}
\begin{equation}\label{eq:Lui}
    \mathcal{L}_{u_i} ^{\rm T} = l_{u_i} ^{\rm T} (y_i) \Delta \tau - \beta_i \Delta \tau + C_{u_i} ^{\rm T} (y_i) \mu_i \Delta \tau + g_{u_i} ^{\rm T} (y_i) \nu_i \Delta \tau = 0,
\end{equation}
and
\begin{equation}\label{eq:duality}
    {\rm diag} (s_i) \nu_i = \epsilon 1, 
\end{equation}
where $\epsilon 1$ denotes an $m_g$-dimensional vector, all of whose elements are $\epsilon$. 
The last equation (\ref{eq:duality}) denotes the complementarity conditions between the slack variable $s_i$ and the Lagrange multiplier $\nu_i$.

\section{Solution Method of Optimal Control Problem}\label{section:solutionmethod}
\subsection{Linearization for Newton's method}
The above KKT conditions are linearized for Newton's method. We apply Gauss-Newton Hessian approximation because the inverse dynamics constraint (\ref{eq:id}) is complicated enough to make it impractical to compute the second-order partial derivatives of (\ref{eq:id}). 
Accordingly, the Hessian of the Lagrangian at the terminal stage is approximated by the Gauss-Newton-approximated Hessian of the terminal cost and that at the intermediate stage by the Gauss-Newton-approximated Hessian of the stage cost.
For example, when the terminal cost and the stage cost take a quadratic form, we have $\mathcal{L}_{q_N q_N} \simeq \varphi_{q_N q_N}$ and $\mathcal{L}_{q_i q_i} \simeq l_{q_i q_i}$.
With the approximated Hessian, (\ref{eq:LqvN}) is linearized into a linear equation with respect to the Newton directions $\Delta \lambda_N, \Delta \gamma_N, \Delta q_N, \Delta v_N \in \mathbb{R}^n$ as 
\begin{align}\label{eq:LqvNLinearized}
    & \begin{bmatrix}
        \mathcal{L}_{q_N} ^{\rm T} \\
        \mathcal{L}_{v_N} ^{\rm T} 
    \end{bmatrix}
    + 
    \begin{bmatrix}
        \mathcal{L}_{q_N q_N} & \mathcal{L}_{q_N v_N} \\
        \mathcal{L}_{v_N q_N} & \mathcal{L}_{v_N v_N}
    \end{bmatrix}
    \begin{bmatrix}
        \Delta q_N \\
        \Delta v_N
    \end{bmatrix} \notag \\
    & + \begin{bmatrix}
        \delta_{q_{N}} ^{\rm T} (q_{N-1}, q_{N}) \Delta \lambda_{N} \\
        - \Delta \gamma_{N}
    \end{bmatrix} = 0,
\end{align}
As well, (\ref{eq:Lqvi})--(\ref{eq:Lui}) are linearized into linear equations with respect to the Newton directions $\Delta \lambda_{i+1}, \allowbreak \Delta \gamma_{i+1}, \allowbreak \Delta q_i , \allowbreak \Delta v_i, \allowbreak \Delta a_i, \allowbreak \Delta u_i, \Delta \beta_i \in \mathbb{R}^n$, $\Delta f_i \in \mathbb{R}^{n_f}$, and $\Delta \mu_i \allowbreak \in \mathbb{R}^{m_c}$.
Note that the directions related to the PDIPM, $\Delta s_i, \Delta \nu_i \in \mathbb{R}^{m_g}$, are eliminated explicitly from the linear equations by adding certain terms related to the logarithmic barrier functions to the Hessians (e.g., ${\mathcal{L}}_{q_i q_i}$) and residuals of the KKT conditions (e.g., ${\mathcal{L}}_{q_i}$) \cite{bib:ipopt, bib:nocedal}. 
By denoting such modified Hessians as $\bar{\mathcal{L}}_{q_i q_i}$ and KKT residuals as $\bar{\mathcal{L}}_{q_i}$, we obtain
\begin{align}\label{eq:LqviLinearized}
    & \begin{bmatrix}
        \bar{\mathcal{L}}_{q_i} ^{\rm T} \\
        \bar{\mathcal{L}}_{v_i} ^{\rm T} 
    \end{bmatrix}
    + \begin{bmatrix}
        \bar{\mathcal{L}}_{q_i, y_i} \\
        \bar{\mathcal{L}}_{v_i, y_i} 
    \end{bmatrix} 
    \Delta y_i
    + \begin{bmatrix}
        \delta_{q_{i}} ^{\rm T} (q_{i}, q_{i+1}) & O_{n \times n} \\
        I_n \Delta \tau & I_n 
    \end{bmatrix}
    \begin{bmatrix}
        \Delta \lambda_{i+1} \\
        \Delta \gamma_{i+1}
    \end{bmatrix} \notag \\
    & + \begin{bmatrix}
        {\rm ID}_{q_i} ^{\rm T} (\tilde{y}_i) \Delta \tau \\
        {\rm ID}_{v_i} ^{\rm T} (\tilde{y}_i) \Delta \tau
    \end{bmatrix}
    \Delta \beta_i  
    + \begin{bmatrix}
        {C}_{q_i} ^{\rm T} ({y}_i) \Delta \tau \\
        {C}_{v_i} ^{\rm T} ({y}_i) \Delta \tau
    \end{bmatrix}
    \Delta \mu_i \notag \\
    & + \begin{bmatrix}
        \delta_{q_{i}} ^{\rm T} (q_{i-1}, q_i) \Delta \lambda_{i} \\
        - \Delta \gamma_{i}
    \end{bmatrix} = 0,
\end{align}
\begin{align}\label{eq:LafiLinearized}
    & \begin{bmatrix}
        \bar{\mathcal{L}}_{a_i} ^{\rm T} \\
        \bar{\mathcal{L}}_{f_i} ^{\rm T} 
    \end{bmatrix}
    + \begin{bmatrix}
        \bar{\mathcal{L}}_{a_i y_i} \\
        \bar{\mathcal{L}}_{f_i y_i} 
    \end{bmatrix}
    \Delta y_i
    + \begin{bmatrix}
        O_{n \times 1} \\
        \Delta \gamma_{i+1} \Delta \tau
    \end{bmatrix} 
    \notag \\
    & 
    + \begin{bmatrix}
        {\rm ID}_{a_i} ^{\rm T} (\tilde{y}_i) \Delta \tau \\
        {\rm ID}_{f_i} ^{\rm T} (\tilde{y}_i) \Delta \tau
    \end{bmatrix}
    \Delta \beta_i 
    + \begin{bmatrix}
        {C}_{a_i} ^{\rm T} ({y}_i) \Delta \tau \\
        {C}_{f_i} ^{\rm T} ({y}_i) \Delta \tau
    \end{bmatrix}
    \Delta \mu_i = 0,
\end{align}
and
\begin{equation}\label{eq:LuiLinearized}
    \bar{\mathcal{L}}_{u_i} ^{\rm T} + \bar{\mathcal{L}}_{u_i y_i} \Delta y_i - \Delta \beta_i \Delta \tau + C_{u_i} ^{\rm T} (y_i) \Delta \mu_i \Delta \tau = 0.
\end{equation}
The constraints, (\ref{eq:id})--(\ref{eq:equalityconstraints}) and (\ref{eq:slack}), are linearized as
\begin{equation}\label{eq:initialstateLinearized}
    \delta(\bar{q}, q_0) + \delta_{q_0} (\bar{q}, q_0) \Delta q_0 = 0, \; \bar{v} - v_0 - \Delta v_0 = 0, 
\end{equation}
\begin{align}\label{eq:forwardEulerLinearized}
    & \begin{bmatrix}
        \delta (q_i, q_{i+1}) + v_i \Delta \tau \\
        v_i - v_{i+1} + a_i \Delta \tau
    \end{bmatrix} 
    + \begin{bmatrix}
        \delta_{q_i} (q_i, q_{i+1}) & I_n \Delta \tau \\
        O_{n \times n} & I_n 
    \end{bmatrix} 
    \begin{bmatrix}
        \Delta {q_i}  \\
        \Delta {v_i}
    \end{bmatrix} 
    \notag \\
    & \;\; 
    + \begin{bmatrix}
        O_{n \times 1} \\
        \Delta a_i \Delta \tau
    \end{bmatrix} 
    + \begin{bmatrix}
        \delta_{q_{i+1}} (q_i, q_{i+1}) & O_{n \times n} \\
        O_{n \times n} & -I_n 
    \end{bmatrix} 
    \begin{bmatrix}
        \Delta {q_{i+1}}  \\
        \Delta {v_{i+1}}
    \end{bmatrix}
    = 0, 
\end{align}
\begin{equation}\label{eq:idLinearized}
    {\rm ID}_{\tilde{y}_i} (\tilde{y}) \Delta \tilde{y}_i
    - \Delta u_i + {\rm ID} (\tilde{y}) - u = 0, 
\end{equation}
and 
\begin{equation}\label{eq:equalityConstraintLinearized}
    {C}_{\tilde{y}_i} (y_i) \Delta \tilde{y}_i
    + {C}_u (y_i) \Delta u_i 
    + C (y_i) = 0.
\end{equation}
Each Newton iteration consists of solving a linear equation that finds Newton directions satisfying (\ref{eq:LqvNLinearized})--(\ref{eq:equalityConstraintLinearized}).

\subsection{Condensing inverse dynamics}
Next, we condense the inverse dynamics; i.e., we eliminate $\Delta u_i$ and $\Delta \beta_i$ from the linear equations (\ref{eq:LqvNLinearized})--(\ref{eq:equalityConstraintLinearized}).
By substituting the expression of $\Delta u_i$ and $\Delta \beta_i$ with respect to other Newton directions (\ref{eq:idLinearized}) and (\ref{eq:LuiLinearized}) into (\ref{eq:LqviLinearized})--(\ref{eq:LuiLinearized}) and (\ref{eq:equalityConstraintLinearized}), we can obtain a condensed linear equation. 
For notational simplicity, we introduce the condensed Hessian,
\begin{align}
    \tilde{\mathcal{L}}_{z_i w_i} := \; & \bar{\mathcal{L}}_{z_i w_i} + {\rm ID}_{z_i} ^{\rm T} (\tilde{y}_i) \bar{\mathcal{L}}_{u_i u_i} {\rm ID}_{w_i} (\tilde{y}_i) \notag \\
    & + {\rm ID}_{z_i} ^{\rm T} (\tilde{y}) \bar{\mathcal{L}}_{u_i w_i} 
    + \bar{\mathcal{L}}_{z_i u_i} {\rm ID}_{w_i} (\tilde{y}) 
\end{align}
for $z_i, w_i \in \left\{ q_i, v_i, a_i, f_i \right\}$ and the condensed KKT residual, 
\begin{align}
    {\tilde{\mathcal{L}}_{z_i}}^{\rm T} := \; & {\bar{\mathcal{L}}_{z_i}}^{\rm T} + {\rm ID}_{z_i} ^{\rm T} (\tilde{y}) {\bar{\mathcal{L}}_{u_i} } ^{\rm T} \notag \\
    & + (\bar{\mathcal{L}}_{z_i, u_i} + {\rm ID}_{z_i} ^{\rm T} (\tilde{y}) \bar{\mathcal{L}}_{u_i u_i}) ({\rm ID} (\tilde{y}) - u_i )
\end{align}
for $z_i \in \left\{ q_i, v_i, a_i, f_i \right\}$. We also introduce the condensed Jacobian of the equality constraint,
\begin{align}
    \tilde{C}_{z_i} := {C}_{z_i} + C_{u_i} {\rm ID}_{z_i}  (q_i, v_i, a_i, f_i) ,
\end{align}
for $z_i \in \left\{ q_i, v_i, a_i, f_i \right\}$ and the condensed residual of the equality constraint, 
\begin{align}
    \tilde{C} := {C} + C_{u_i} ({\rm ID}  (q_i, v_i, a_i, f_i) - u_i).
\end{align}
Then, $\Delta u_i$ and $\Delta \beta_i$ are eliminated from (\ref{eq:LqviLinearized})--(\ref{eq:LuiLinearized}) and (\ref{eq:equalityConstraintLinearized}):
\begin{align}\label{eq:LqviCondensed}
    & \begin{bmatrix}
        \tilde{\mathcal{L}}_{q_i} ^{\rm T} \\
        \tilde{\mathcal{L}}_{v_i} ^{\rm T} 
    \end{bmatrix}
    + \begin{bmatrix}
        \tilde{\mathcal{L}}_{q_i y_i} \\
        \tilde{\mathcal{L}}_{v_i y_i}
    \end{bmatrix} 
    \Delta y_i 
    + \begin{bmatrix}
        \delta_{q_{i}} ^{\rm T} (q_{i}, q_{i+1}) & O_{n \times n} \\
        I_n \Delta \tau & I_n 
    \end{bmatrix}
    \begin{bmatrix}
        \Delta \lambda_{i+1} \\
        \Delta \gamma_{i+1}
    \end{bmatrix}
    \notag \\
    & + \begin{bmatrix}
        {C}_{q_i} ^{\rm T} ({y}_i) \\
        {C}_{v_i} ^{\rm T} ({y}_i)
    \end{bmatrix}
    \Delta \mu_i + \begin{bmatrix}
        \delta_{q_{i}} ^{\rm T} (q_{i-1}, q_i) \Delta \lambda_{i} \\
        - \Delta \gamma_{i}
    \end{bmatrix} = 0,
\end{align}
\begin{equation}\label{eq:LafiCondensed}
    \begin{bmatrix}
        \tilde{\mathcal{L}}_{a_i} ^{\rm T} \\
        \tilde{\mathcal{L}}_{f_i} ^{\rm T} 
    \end{bmatrix}
    + \begin{bmatrix}
        \tilde{\mathcal{L}}_{a_i y_i} \\
        \tilde{\mathcal{L}}_{f_i y_i} 
    \end{bmatrix}
    \Delta y_i
    + \begin{bmatrix}
        O_{n \times 1} \\
        \Delta \gamma_{i+1} \Delta \tau
    \end{bmatrix}
    + \begin{bmatrix}
        {C}_{a_i} ^{\rm T} ({y}_i) \\
        {C}_{f_i} ^{\rm T} ({y}_i)
    \end{bmatrix}
    \Delta \mu_i = 0,
\end{equation}
and
\begin{equation}\label{eq:CCondensed}
    \tilde{C}_{\tilde{y}_i} (y_i) \Delta \tilde{y}_i + \tilde{C} (y_i) = 0.
\end{equation}
After condensing, the linear equation is reduced to find the Newton directions $\Delta \lambda_i, \Delta \gamma_i, \Delta q_i , \Delta v_i, \Delta a_i, \Delta f_i, \Delta \mu_i$ satisfying (\ref{eq:initialstateLinearized}), (\ref{eq:forwardEulerLinearized}), (\ref{eq:LqvNLinearized}), (\ref{eq:LqviCondensed})--(\ref{eq:CCondensed}). 

If the rigid-body system is fully actuated and there are no contacts, the size of the condensed linear equation is the same as that of the multiple-shooting method based on forward dynamics. 
If it is underactuated and there are no contacts, the size of the condensed linear equation at each time stage is increased from that of the multiple-shooting method based on forward dynamics by only twice the number of the passive joints, since we assume the system is fully actuated and add an equality constraint to zero the control input torques corresponding to the passive joints. 
If there are contacts, the size of the linear equation also increases depending on the way to treat the contacts.

\subsection{Algorithm}\label{subsec:algorithm}
Algorithm \ref{alg:idocp} is the pseudocode of a single Newton iteration of the proposed method. 
The first step (lines 1--3) forms the linear equation consisting of (\ref{eq:initialstateLinearized}), (\ref{eq:forwardEulerLinearized}), (\ref{eq:LqvNLinearized}), and (\ref{eq:LqviCondensed})--(\ref{eq:CCondensed}) by computing the condensed KKT residuals and Hessians, such as $\tilde{\mathcal{L}}_{q_i}$, $\tilde{\mathcal{L}}_{q_i q_i}$, $\tilde{C}$, and $\tilde{C}_q$. 
This step is fully parallelizable into each time stage. The second step (line 4--6) computes the directions $\Delta \lambda_i$, $\Delta \gamma_i$, $\Delta \tilde{y}_i$, $\Delta \mu_i$ by solving the linear equation consisting of (\ref{eq:initialstateLinearized}), (\ref{eq:forwardEulerLinearized}), (\ref{eq:LqvNLinearized}), and (\ref{eq:LqviCondensed})--(\ref{eq:CCondensed}), i.e., inverting the condensed Hessian. 
In this step, we can utilize various efficient methods tailored for the OCP \cite{bib:MPCSurvey}, e.g., Riccati recursion \cite{bib:Riccati} and a highly parallelizable method \cite{bib:parallelNewton}. 
The third step (lines 7--9) computes the condensed directions, i.e., $\Delta u_i$ and $\Delta \beta_i$ from (\ref{eq:idLinearized}) and (\ref{eq:LuiLinearized}), and $\Delta s_i$ and $\Delta \nu_i$ according to the PDIPM \cite{bib:ipopt, bib:nocedal}.
The fourth step (line 10) determines the step size, e.g., by using the fraction-to-boundary rule \cite{bib:nocedal, bib:ipopt}. 
Finally, all variables are updated according to the step size (lines 11--16). 

\begin{figure}[!t]
\begin{algorithm}[H]
\caption{Single Newton iteration of primal-dual interior point method with condensing of inverse dynamics}
\label{alg:idocp}
\begin{algorithmic}[1]
    \Require Initial state $\bar{q}, \bar{v}$
    \Ensure $\lambda_0, ..., \lambda_N$, $\gamma_0, ..., \gamma_N$, $y_0, ..., y_{N-1}$, $q_N$, $v_N$, $\mu_0, ..., \mu_{N-1}$, $\beta_0, ..., \beta_{N-1}$, $s_0, ..., s_{N-1}$, $\nu_0, ..., \nu_{N-1}$
    \For{$i=0,\cdots,N$} {\bf in parallel}
        \State Compute the condensed Hessian and KKT residual.
    \EndFor
    \For{$i=0,\cdots,N$} {\bf in parallel or serial}
        \State Compute the Newton directions $\Delta \lambda_i$, $\Delta \gamma_i$, $\Delta \tilde{y}_i$, $\Delta \mu_i$ by solving the linear equation (\ref{eq:initialstateLinearized}), (\ref{eq:forwardEulerLinearized}), (\ref{eq:LqvNLinearized}),  (\ref{eq:LqviCondensed})--(\ref{eq:CCondensed}).
    \EndFor
    \For{$i=0,\cdots,N-1$} {\bf in parallel}
        \State Compute the condensed Newton directions $\Delta u_i$ and $\Delta \beta_i$ from (\ref{eq:idLinearized}), (\ref{eq:LuiLinearized}), 
        and $\Delta s_i$ and $\Delta \nu_i$ according to \cite{bib:ipopt}.
    \EndFor
    \State Determine the step size $\alpha \in (0, 1]$.
    \For{$i=0,\cdots,N$} {\bf in parallel}
        \State Update all variables $\lambda_i$, $\gamma_i$, $y_i$, $\mu_i$ $\beta_i$, $s_i$, $\nu_i$ by 
        \State $\lambda_i \leftarrow \lambda_i + \alpha \Delta \lambda_i$, $\gamma_i \leftarrow \gamma_i + \alpha \Delta \gamma_i$, $y_i \leftarrow y_i + \alpha \Delta y_i$,
        \State $\mu_i \leftarrow \mu_i + \alpha \Delta \mu_i$, $\beta_i \leftarrow \beta_i + \alpha \Delta \beta_i$, 
        \State $s_i \leftarrow s_i + \alpha \Delta s_i$, $\nu_i \leftarrow \nu_i + \alpha \Delta \nu_i$
    \EndFor
\end{algorithmic}
\end{algorithm}
\vspace{-5mm}
\end{figure}

The main advantage of the proposed method appears in the first step (lines 1--3). It calls RNEA and the partial derivatives of RNEA $N-1$ times, whereas the forward dynamics-based method computes ABA (or inverse of the joint-inertia matrix) $N-1$ times and the partial derivatives of ABA $N-1$ times. Moreover, the proposed method can reduce the computational cost in the second step (lines 4--6), as our formulation leads to sparsity in the partial derivatives of the state equation. 
We will explain this point below with the Riccati recursion for computing the Newton directions.
Riccati recursion \cite{bib:Riccati} performs a recursive block elimination to compute the Newton directions with the forward Euler discretization (\ref{eq:forwardEuler}). 
In Riccati recursion, we regard $[\Delta q_i ^{\rm T} \;\; \Delta v_i ^{\rm T} ]^{\rm T}$ as the state variable and $[\Delta a_i ^{\rm T} \;\; \Delta f_i ^{\rm T} ]^{\rm T}$ as the control input of a linear quadratic regulator subproblem.
We also introduce matrices, 
\begin{equation*}\label{eq:RiccatiFactorization}
    A_i := \begin{bmatrix}
        \delta_{q_i} (q_{i}, q_{i+1}) & I_n \Delta \tau  \\
        O_{n \times n} & I_n
    \end{bmatrix}, \;\;
    B_i := \begin{bmatrix}
        O_{n \times n} & O_{n \times n_f}  \\
        \Delta \tau I_n & O_{n \times n_f} 
    \end{bmatrix},
\end{equation*}
where $A_i$ is composed of the partial derivatives of the state equation with respect to $q$ and $v$, and $B_i$ is composed of the partial derivatives with respect to $a$ and $f$. Riccati recursion serially computes $A_i^{\rm T} P_{i+1} A_i$, $A_i ^{\rm T} P_{i+1} B_i$, $B_i ^{\rm T} P_{i+1} B_i$ for a sequence of matrices $P_{i}$ that represent the sensitivities of $\Delta \lambda_i$ and $\Delta \gamma_i$ with respect to $\Delta q_i$ and $\Delta v_i$. Thanks to the sparsity of $A_i$ and $B_i$ in the proposed formulation, several matrix products reduce to sums of matrices, and this in turn reduces the computational cost, especially the cost of the serial computational part.


We have developed an open-source C++ software framework for the OCP for rigid-body systems, \texttt{idocp} \cite{bib:idocpWeb}, that utilizes the above Riccati recursion and ParNMPC, a highly parallelizable method \cite{bib:parallelNewton}. 
\texttt{idocp} uses Eigen\footnote{http://eigen.tuxfamily.org} for linear algebra and \texttt{pinocchio} \cite{bib:pinocchio}, an efficient C++ library for rigid-body dynamics algorithms, to compute RNEA and its partial derivatives. 
It also employs the parallel Newton-type method \cite{bib:parallelNewton} with backward Euler discretization method.

\section{Numerical Experiments}\label{section:expriment}
\subsection{Experimental settings}
To evaluate the computational efficiency and numerical robustness of the proposed solution method of the OCP, we compared our solver \texttt{idocp} with \texttt{crocoddyl} \cite{bib:crocoddyl}, a highly efficient C++ implementation of DDP/iLQR for rigid-body systems, and Ipopt \cite{bib:ipopt}, an off-the-shelf nonlinear optimization solver.
Both \texttt{idocp} and \texttt{crocoddyl} utilize \texttt{pinocchio} for rigid-body dynamics algorithms, Eigen for linear algebra, and OpenMP \cite{bib:OpenMP} for parallel computing. 
Ipopt was customized to solve the OCP for rigid-body systems based on forward dynamics and the multiple-shooting method with parallel computing by using \texttt{pinocchio} and OpenMP.
As the algorithm of \texttt{idocp}, we used Riccati recursion to compute the Newton directions (we will refer to this method as inverse dynamics-based Riccati recursion (IDRR) hereafter). Furthermore, we utilized two algorithms from \texttt{crocoddyl}: the standard DDP with a Gauss-Newton Hessian approximation, which is identical to iLQR, and feasible-prone DDP (FDDP) \cite{bib:crocoddyl}, a kind of iLQR that improves numerical robustness by modifying the backward pass of iLQR in a multiple-shooting fashion \cite{bib:multipleshootingiLQR} and utilizes a line-search method based on the Goldstein condition \cite{bib:nocedal}. We will refer to the former as iLQR and the latter as FDDP. 
In Ipopt, we used the Broyden–Fletcher–Goldfarb–Shanno method for the Hessian approximation and Harwell Subroutine Library MA57 to solve the linear subproblems. 
All experiments were conducted on a laptop with quad-core CPU Intel Core i7-10510U @1.8GHz.

\subsection{Computational time}\label{subsub:CPUTime}

\begin{figure}[tb]
\centering
\includegraphics[scale=0.62]{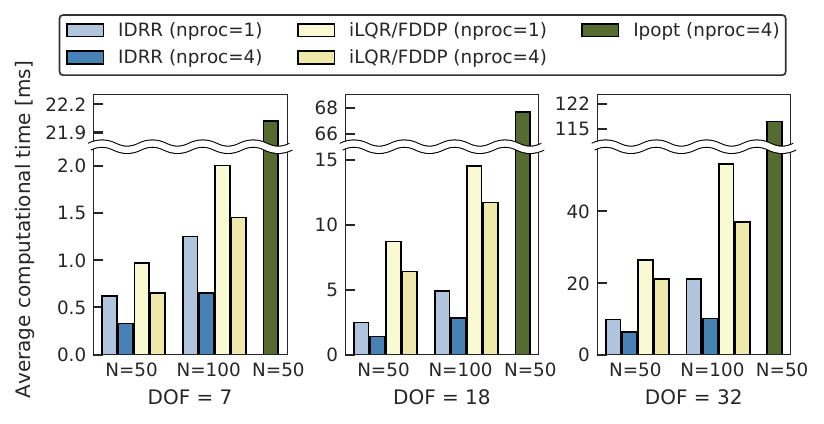}
\caption{Average computational time per iteration ${\rm [ms]}$ of the proposed method (${\rm IDRR}$), ${\rm iLQR}$/${\rm FDDP}$, and Ipopt, for various degrees of freedom (DOF $\in \left\{ 7, 14, 32 \right\}$), numbers of time steps ($N \in \left\{ 50, 100 \right\}$), and numbers of threads (${\rm nproc} \in \left\{ 1, 4 \right\}$)}
\label{fig:time}
\end{figure}

First, we compared the computational time per iteration of IDRR, iLQR/FDDP, and Ipopt. 
In particular, we computed the average computational time over 10,000 trials for systems with various degrees of freedom (DOF), various numbers of time stages $N$, and various numbers of threads (nproc) with the quadratic terminal cost,
\begin{equation}\label{eq:terminalcost}
    \varphi (q, v) = \frac{1}{2} q_e ^{\rm T} Q_q q_e + \frac{1}{2} v_e ^{\rm T} Q_v v_e , 
\end{equation}
where $Q_q = I_{n}$, $Q_v = I_{n}$, $q_e := q - q_{\rm ref}$, and $v_e := v - v_{\rm ref}$ with $q_{\rm ref}, v_{\rm ref} \in \mathbb{R}^n$, and the quadratic stage cost, 
\begin{equation}\label{eq:stagecost}
    l (q, v, a, u) = \frac{1}{2} q_e ^{\rm T} Q_q q_e 
    + \frac{1}{2} v_e ^{\rm T} Q_v v_e
    + \frac{1}{2} u_e ^{\rm T} Q_u u_e, 
\end{equation}
where $Q_u = 0.001 \times I_{n}$ and $u_e := u - u_{\rm ref}$ with $u_{\rm ref} \in \mathbb{R}^n$.
Fig.\ \ref{fig:time} shows the computational time per iteration [ms] of each solver. Note that the results of iLQR and FDDP are shown as one because they had almost the same computational time. 
Moreover, the results of Ipopt are only for the fastest case, i.e., the case with $N=50$ and nproc = 4, because Ipopt was much slower than the other solvers.
As shown in the figure, for all combinations of $N$ and number of threads (nproc), our IDRR was faster than iLQR/FDDP and Ipopt. This demonstrates that the inverse dynamics-based formulation can reduce the computational cost compared with the forward dynamics-based formulation. We also found that IDRR became faster as the number of the threads increased, whereas iLQR/FDDP did so only moderately. This is because the proposed method reduces the computational burden of the serial calculation in the Riccati recursion, thanks to its multiple-shooting and sparsity structure, as stated in \ref{subsec:algorithm}.

\subsection{Numerical robustness}

\begin{figure}[tb]
    \centering
    \includegraphics[scale=0.64]{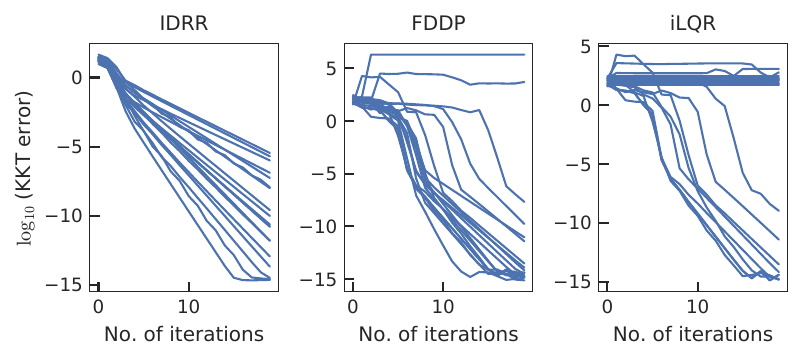}
    \caption{
    $\log_{10}$ scaled KKT errors of the proposed method (IDRR), FDDP, and iLQR over 20 random trials. 
    }
    \label{fig:convergence}
\end{figure}

Second, we investigated the numerical robustness, i.e., the convergence, of the proposed method through the OCP for KUKA iiwa14, a 7-DOF manipulator.
We set the length of the horizon to 1 s and divided it into $N=50$ steps, i.e., $\Delta \tau = 0.02 {\rm [s]}$. 
The objective of the OCP is to make the configuration $q$ converge to $q_{\rm ref}$ and the velocity $v$ converge to zero given a random initial state $\bar{q}$ and $\bar{v}$. We used the quadratic terminal cost (\ref{eq:terminalcost}) and stage cost (\ref{eq:stagecost}) and set $q_{\rm ref}$ to $[0 \;\; \pi / 2 \;\; 0 \;\; \pi / 2 \;\; 0 \;\; \pi / 2 \;\; 0]^{\rm T}$, $v_{\rm ref}$ to zero, and $u_{\rm ref}$ to the gravity compensation torques at $q_{\rm ref}$.
We did not impose any inequality constraints in this experiment. 
iLQR and FDDP used a line search and did not use regularization. IDRR did not use either, while Ipopt used both of them.
We initialized $q_i$ and $v_i$ in the solution of IDRR, FDDP, and Ipopt by $\bar{q}$ and $\bar{v}$. We randomly selected each element of $\bar{q}$ from $[-1, 1]$ and each element of $\bar{v}$ from $[-10, 10]$. We ran 20 trials for each solver and computed the $l^2$-norm of the residual of the KKT conditions, which we will refer to as the KKT error hereafter, and the total cost to be minimized. Note that the KKT conditions of IDRR are given by (\ref{eq:id})--(\ref{eq:equalityconstraints}), (\ref{eq:slack}), and (\ref{eq:LqvN})--(\ref{eq:duality}), while those of iLQR/FDDP and Ipopt are composed of different equations based on forward dynamics and/or single-shooting. Figure \ref{fig:convergence} shows the $\log_{10}$-scaled KKT errors of IDRR, FDDP, and iLQR. 
Results for Ipopt are not shown because it converged rather slowly (over 1000 iterations were needed to reduce the KKT error to under $10^{-10}$ in most cases) due to the quasi-Newton Hessian approximation.
Note as well that, as the equations representing the KKT conditions differ between IDRR and FDDP/iLQR, we cannot directly compare their KKT errors. However, from the graph of FDDP, we can clearly see that FDDP diverged in two trials. In addition, from the graph of iLQR, we can see that in many trials the KKT error could not be completely reduced because the iLQR did not produce a descent direction and the step size became zero. On the other hand, IDRR reduced the KKT error at each iteration in all cases, including those in which FDDP and iLQR failed to converge. 
This indicates that the proposed method is more robust than single-shooting methods such as iLQR and FDDP.

\subsection{MPC for floating base systems}
\begin{figure}[tb]
    \centering
    \includegraphics[scale=0.21]{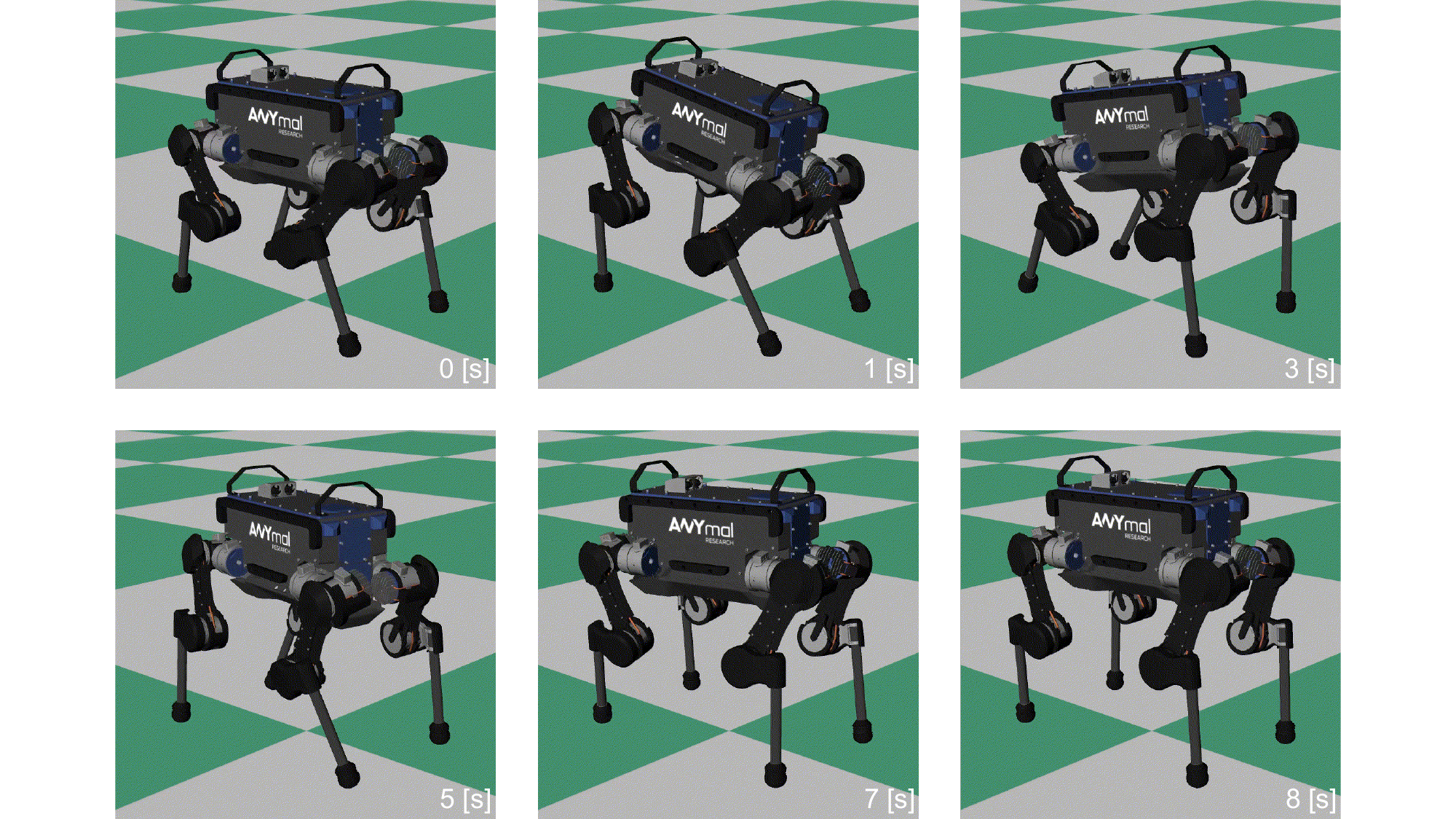}
    \caption{Posture control of ANYmal by whole-body MPC controller.}
    \label{fig:anymal:postures}
\end{figure}

Finally, we investigated the applicability of the proposed method to floating base systems, which are a kind of underactuated system, by simulating MPC implemented by IDRR for the whole-body control of a quadruped ANYmal. To take the rigid contacts into account, we treated the contact forces as the optimization variables and imposed an equality constraint in the form of Baumgarte's stabilization method \cite{bib:baumgarte} for each contact, i.e., a 12-dimensional equality constraint for four contacts. 
We used a quadratic cost function to track the desired configuration. 
We also imposed inequality constraints on the joint angle limits, angular velocity limits, and torque limits and linearized friction cones using the PDIPM.
We set the length of the horizon to 1 s and divided it into $N=20$ equal steps. We set the sampling time to 2.5 ms, and the MPC controller updated the solution once per sampling period. We utilized \texttt{RaiSim} \cite{bib:raisim}, an articulated-body simulator with contacts. Figure \ref{fig:anymal:postures} shows the various postures of ANYmal controlled by MPC. As can be seen, the proposed method was able to control the underactuated system. Each control update took around 2.1 ms with four threads on the same laptop as in \ref{subsub:CPUTime} and achieved real-time MPC. 

\section{Conclusions}\label{section:conclu}
We proposed an efficient method of solving the OCP for rigid-body systems on the basis of inverse dynamics and the multiple-shooting method. In this method, we regard all variables, including the state, acceleration, and control input torques, as optimization variables and treat the inverse dynamics as an equality constraint. We eliminated the update in the control input torques from the linear equation of Newton's method by applying condensing for inverse dynamics. The size of the resultant linear equation was the same as that of the multiple-shooting method based on forward dynamics except for the variables related to the passive joints and contacts. The proposed method reduces the computational cost of the dynamics and their sensitivities by utilizing RNEA and its partial derivatives. In addition, it increases the sparsity of the Hessian of the KKT conditions, which further reduces the computational cost, e.g., of Riccati recursion. Numerical experiments showed that the proposed method is more than twice as fast as iLQR and a DDP variant based on forward dynamics. They also showed that it is more numerically robust than these conventional methods.

Our future work will include incorporating changes in contact status by means of the hybrid optimal control approach \cite{bib:hybridMPC} or the contact-implicit approach \cite{bib:todorov:mujoco, bib:implicitContact} with an eye toward developing applications for manipulation and locomotion. 
Both methods can be combined with the proposed method because the contact forces are included in the optimization variables.

\addtolength{\textheight}{-9cm}   







\bibliographystyle{IEEEtran}
\bibliography{IEEEabrv, ieee}


\end{document}